\documentclass[12pt]{article}

\def\update{April 21, 2008}

\usepackage{
amssymb,
amsthm}

\usepackage{amsmath}
\usepackage{amsfonts}%para \mathbb{} Z, R, C y N

 \usepackage[colorlinks,linkcolor=black]{hyperref}

\newtheorem{definition}{\textbf{Definition}}

\newtheorem{theorem}{\textbf{Theorem}}
\newtheorem{corollary}{\textbf{Corollary}}
\newtheorem{remark}{\textbf{Remark}}
\newtheorem{conjecture}{\textbf{Conjecture}}
\newtheorem{lemma}{\textbf{Lemma}}

\def\Z {\mathbb{Z}}
\def\Q {\mathbb{Q}}
\def\QQ {\overline{\Q}}

\def\C {\mathbb{C}}

\def\fle {\longrightarrow }

\def\trdeg{\mathrm{trdeg}}

\begin{document}

\noindent
AWS 2008 Project
\hfill
{updated: \it \update}

\begin{center}

\medskip

{
\LARGE\mathversion{bold}
Some consequences of Schanuel's Conjecture
}
\medskip

 \large

\medskip

{\it by }
\medskip

Chuangxun Cheng, Brian Dietel,

Mathilde Herblot,   Jingjing Huang, Holly Krieger,

 Diego Marques,  Jonathan Mason, Martin Mereb

 and S. Robert Wilson.

\end{center}

During the Arizona Winter School 2008 (held in Tucson, AZ) we
worked on the following problems:

a) (Expanding a remark by S. Lang  \cite{La1}). Define $E_0 = \QQ$.
Inductively, for $n \geq 1$, define $E_n$ as the algebraic closure
of the field generated over $E_{n-1}$ by the numbers
$\exp(x)=e^x$, where $x$ ranges over $E_{n-1}$. Let $E$ be the
union of $E_n$, $n \geq 0$. Show that Schanuel's Conjecture
implies that the numbers $\pi, \log \pi, \log \log \pi, \log \log
\log \pi, \ldots $ are algebraically independent over $E$.

 b) Try to get a (conjectural) generalization
involving the field $L$ defined as follows. Define $L_0 = \QQ$.
Inductively, for $n \geq 1$, define $L_n$ as the algebraic closure
of the field generated over $L_{n-1}$ by the numbers $y$, where
$y$ ranges over the set of complex numbers such that $e^y\in
L_{n-1}$. Let $L$ be the union of $L_n$, $n \geq 0$.

We were able to prove the more general result:

\begin{theorem} Schanuel's Conjecture implies $E$ and $L$ are linearly
disjoint over $\overline{\Q}$.

\end{theorem}

And deduced from it the following ones:

\begin{enumerate}
    \item $\pi \not\in E$ and $e \not\in L$.

    \item $\pi, \log \pi, \log\log\pi, \ldots$ are
    algebraically independent over $E$.

    \item $e, e^e, e^{e^e}, \ldots$ are algebraically independent over
    $L$.

    \item $E\cap L=\overline{\Q}$.
\end{enumerate}

\newpage %pagebreak?

Remember:

\begin{conjecture} [Schanuel] Let $ x_1, \ldots , x_n$ be $\Q$-linearly
independent complex numbers. Then the transcendence degree over
$\Q$ of the field
$$
\Q( x_1 , \ldots , x_n , e^{x_1} , \ldots ,
e^{x_n} )
$$
 is at least $n$.
\end{conjecture}

\begin{definition} Let $F/K$ be a field extension and $F_1,F_2\subseteq F$ two subextensions. We say they are \emph{linearly disjoint}
over $K$ when the following holds:

 $\{x_1, \ldots, x_n\}\subseteq F_1$ linearly
independent over $K$ $\Rightarrow$  $\{x_1, \ldots, x_n\}$
linearly independent over $F_2$.

We say they are \emph{free} (or \emph{algebraically disjoint})
over $K$ when:

 $\{x_1, \ldots, x_n\}\subseteq F_1$ algebraically
independent over $K$ $\Rightarrow$ $\{x_1, \ldots, x_n\}$
algebraically independent over $F_2$.
\end{definition}

\begin{remark}
Linear disjointness is equivalent to the multiplication
map
$$
F_1\otimes_K F_2 \fle F
$$
 being injective. Therefore this is
a symmetric condition in $F_1$ and $F_2$.
\end{remark}

\begin{remark}  Algebraic disjointness is equivalent to the existence
of transcendence basis $B_1$, $B_2$ of the extensions $F_1/K$ and
$F_2/K$ (respectively) such that $B_1\cup B_2$ is algebraically
independent over $K$. Therefore this one is also a symmetric
condition in $F_1$ and $F_2$.
\end{remark}

\begin{remark}  For a set $S\subseteq F_1$ to be algebraically independent
over $K$ means all its monomials being linearly independent over
$K$. Thus linearly disjointness \emph{implies} freeness, and in
general the converse is not true (although we are going to use a
partial converse to this fact, proved in \cite{La3}).
\end{remark}

\begin{remark}  If $F_1,F_2$ are linearly disjoint over $K$ then we must
have $F_1\cap F_2 = K$, since $k\in F_1\cap F_2$ will be
$F_2$-linearly dependent together with $1$ whence, they should
also be $K$-linearly dependent.
\end{remark}

Before going to the proof of the Theorem, we need a couple of
technical lemmas involving a key construction.

\begin{lemma} We have $E_n=\overline{\Q(\exp(E_{n-1}))}$\footnote{$\overline{F}$
meaning the algebraic closure of the field $F$.}.
\end{lemma}

\begin{proof}
 With induction in $n$, the base case follows by definition
since
$$
E_1=\overline{E_0(\exp(E_0))}=\overline{\QQ(\exp(E_0))}=\overline{\Q(\exp(E_0))},
$$

and $E_0=\QQ$.

In general
\begin{eqnarray}
E_n & = &
\overline{E_{n-1}(\exp(E_{n-1}))}\nonumber\\
E_n & = &
\overline{\overline{\Q(\exp(E_{n-2}))}(\exp(E_{n-1}))} \nonumber\\
E_n & = &
\overline{\Q(\exp(E_{n-2}))(\exp(E_{n-1}))}\nonumber\\
E_n & = & \overline{\Q(\exp(E_{n-1}))},
\end{eqnarray}

since $E_{n-2}\subseteq E_{n-1}$. %end of proof
\end{proof}

\begin{lemma} For every $x\in E_n$ there is a finite set
$A_{n-1}\subseteq E_{n-1}$ such that $x$ is algebraic over
$\Q(\exp(A_{n-1}))$ (or equivalently,
$x\in\overline{\Q(\exp(A_{n-1}))}$ ).
\end{lemma}

\begin{proof}
We have  $x\in E_n=\overline{\Q(\exp(E_{n-1}))}$ which means it is a
root of a nontrivial polynomial with coefficients in
$\Q(\exp(E_{n-1}))$. Each coefficient involves only finitely many
exponentials of elements in $E_{n-1}$. Therefore, taking $A_{n-1}$
the union of those exponents will work. %end of proof
\end{proof}

\begin{lemma} [the Key Lemma] For every $x\in E_n$ there is a finite set
$A\subseteq E_{n-1}$ such that $x\in\overline{\Q(\exp(A))}$ and
$A$ is also algebraic over $\Q(\exp(A))$.
\end{lemma}

\begin{proof}
 Start with $A_{n-1}$ as in the previous lemma and iterate
the reasoning finding a sequence of subsets
$A_{n-1},A_{n-2},A_{n-3},\ldots,A_0$ as follows:

\begin{itemize}
  \item Since $A_{n-1}\subseteq E_{n-1}$ is finite, it follows that  $A_{n-1}$ is algebraic over $\Q(\exp(A_{n-2}))$ for some finite $A_{n-2}\subseteq
  E_{n-2}$.

  \item Next $A_{n-2}$ is algebraic over $\Q(\exp(A_{n-3}))$ for some finite $A_{n-3}\subseteq E_{n-3}$.
\ldots
  \item Finally $A_{1}$ is algebraic over
$\Q(\exp(A_{0}))$ for some finite $A_{0}\subseteq
E_{0}=\overline{\Q}$.
\end{itemize}

Then just take $A=\bigcup_{m\leq n-1} A_m \subseteq E_{n-1}$.
Since $A_{n-1}\subseteq A$ we get $x\in\overline{\Q(\exp(A))}$ and
since each $A_m$ is algebraic over $\Q(\exp(A_{m-1}))$ then it is
so over $\Q(\exp(A))$ and therefore, the whole set $A$ is
algebraic over $\Q(\exp(A))$. %end of proof
\end{proof}

In a similar way we get analogues of these lemmas in the case of
the logarithmic extensions $L_m$. Let us state them for the sake of
preciseness.

\begin{lemma}
We have $L_n=\overline{\Q(\exp^{-1}(L_{n-1}))}$.
\end{lemma}

\begin{lemma} For every $x\in L_n$ there is a finite set $C_n\subseteq\C$
such that $\exp(C_n)\subseteq L_{n-1}$ and that $x$ is algebraic
over $\Q(C_n)$.
\end{lemma}

\begin{lemma} [the Key Lemma] For every $x\in L_n$ there is a finite set
$C\subseteq\C$ with $\exp(C)\subseteq L_{n-1}$ such that $\exp(C)
\cup \{x\}$ is algebraic over $\Q(C)$.
\end{lemma}

The proofs follow the same outline as in the exponential case.
%should we write them down precisely?

Now we are ready to go the proof of the theorem:

Assuming the Schanuel's Conjecture to be true, let us prove $E_m$
and $L_n$ are linearly disjoint for arbitrary $m$ and $n$ (this
will be enough since $E$ is the union of the $E_m$ and $L$ is the
union of the $L_n$).

Proceeding by induction, let us  assume it is true that $E_{m-1}$
and $L_{n}$ are linearly disjoint over $\QQ$.

Suppose  $E_m$
and $L_n$ are not  linearly disjoint. Let us take $\{l_i\}\subseteq L_n$ linearly
independent over $\overline{\Q}$ and $\{e_i\}\subseteq E_m$ such
that $\sum l_i e_i=0$.

By the Key Lemmas:
\begin{itemize}

\item $\exists$ finite $A \subseteq E_{m-1}$ such that
$A\cup\{e_i\}$ algebraic over $\Q(\exp(A))$.

\item $\exists$ finite $C \subseteq L_{n}$ finite such that
$\exp(C)\cup\{l_i\}$ algebraic over $\Q(C)$.

\end{itemize}

Now take $B\subseteq A$ such that $\exp(B)$ is a transcendence
basis of $\Q(\exp(A))$, and take $D\subseteq C$ such that $D$ is a
transcendence basis of $\Q(C)$.

We claim $B \cup D$ is linearly independent over $\Q$.

Consider any linear relation over $\Q$ and by cleaning
denominators if necessary take
$$
\sum_{b\in B}p_b b=\sum_{d\in D}q_d d
$$
with all the $p_b,q_d$ integers.

Since the expression on the left is an element in $E_{m-1}$ and that of the right is an element of $L_n$, and by hypothesis
these two fields were linearly disjoint over $\QQ$, we should have
$E_{m-1}\cap L_n=\QQ$ and both expressions would represent an
element $r\in\QQ$.

But $\sum_{d\in D}q_d d=r$ is an algebraic relation of $D$ with
coefficients in $\QQ$, hence it must be the trivial relation (keep
in mind that $D$ was taken to be algebraically independent over
$\QQ$).

We get at once $r=0=q_d\:  \; \hbox{ for all } \; d\in D$.

Now from $\sum_{b\in B}p_b b=0$ taking exponentials on both sides
we get
$$
\prod_{b\in B}(\exp(b))^{p_b}=1
$$
 which is an algebraic
relation with coefficients in $\Q$ (and hence in $\QQ$) among the
elements of the set $\exp(B)$, taken algebraically independent.
Therefore, the (Laurent) monomial $\prod_{b\in B}(X_b)^{p_b}$
should be the trivial one, i.e. the integers $p_b$ must be all
equal to zero.

Summarizing, we have proven $B\cup D$ is $\Q$-linearly independent.

By Schanuel's Conjecture $\trdeg_\Q\Q(B,D, \exp(B), \exp(D))\geq
|B|+|D|$.

However
$$
\trdeg_\Q\Q(B,D, \exp(B), \exp(D)) = \trdeg_\Q\Q(B,C, \exp(A), \exp(D))
$$
since $\exp(A)$ is algebraic over $\Q(\exp(B))$ and $C$ is
algebraic over $\Q(D)$.

We also have
$$
\trdeg_\Q\Q(B,C, \exp(A), \exp(D))=\trdeg_\Q\Q(C, \exp(A))
$$
because $B\subseteq A$ and the latter was algebraic over
$\Q(\exp(A))$, and similarly $\exp(D)\subseteq\exp(C)\subseteq
\overline{\Q(C)}$.

Finally
$$
\trdeg_\Q\Q(C, \exp(A))=\trdeg_\Q\Q(D, \exp(B))\leq |B|+|D|
$$
since $C$ was algebraic over $\Q(D)$ and $\exp(A)$ was so over
$\Q(\exp(B))$.

From
$$
|B|+|D| \geq \trdeg_\Q\Q(D, \exp(B))\geq |B|+|D|
$$
 we
conclude $\trdeg_\Q\Q(D, \exp(B))= |B|+|D|$ and the set
$\exp(B)\cup D$ turns out to be algebraically independent over
$\Q$, whence over $\QQ$.

Therefore $\Q(\exp(B))$ and $\Q(D)$ are free over $\overline{\Q}$,
and the same is true for $\overline{\Q(\exp(B))}$ and
$\overline{\Q(D)}$.

Since $\overline{\Q}$ is algebraically closed,
$\overline{\Q(\exp(B))}$ and $\overline{\Q(D)}$ are linearly
disjoint over $\overline{\Q}$ (see \cite{La3} Theorem 4.12, page
367).

But the $\{l_i\}$ are algebraic over $\Q(C)$ and the $\{e_i\}$
are algebraic over $\Q(\exp(A))$, which means $\{l_i\}\subseteq
\overline{\Q(D)}$ and $\{e_i\}\subseteq \overline{\Q(\exp(B))}$
giving to us the nontrivial linear relation $\sum l_i e_i=0$.
Contradiction. %end of proof

\begin{corollary}
We have
 $L\cap E = \QQ$.
\end{corollary}

\begin{proof}
It follows directly from the linear disjointness. %end of proof
\end{proof}

\begin{corollary}
We have
 $\pi \not\in E$ and $e \not\in L$.
\end{corollary}

\begin{proof}
We have
 $e=\exp(1)\in E_1\subseteq E$. Since $e\not\in\QQ$ it
cannot be also in $L$.

If $\pi$ were in $E$, $i\pi$ should also be there. But $i\pi\in
L_1\subseteq L$ since it is a logarithm of $-1$. We conclude
$i\pi\not\in E$ because it is not in $\QQ$. %end of proof
\end{proof}

\begin{corollary}
The numbers  $\pi,\log\pi,\log\log\pi,\ldots$ are algebraically
independent over $E$.
\end{corollary}

\begin{proof}
 We are actually going to prove that
$i\pi,\log\pi,\log\log\pi,\ldots$ are algebraically independent
over $E$ (which is an equivalent statement).

Let us write $\log_{[k]}\pi$ for the $k^{th}-$ iterated logarithm of
$\pi$.

Observe that the whole sequence $i\pi,\log\pi,\log\log\pi,\ldots$
lies in $L$.

Since we are assuming $E$ and $L$ linearly independent over $\QQ$,
they are going be free, and it will be enough to prove
$i\pi,\log\pi,\log\log\pi,\ldots$ they are algebraically
independent over $\QQ$, or, which is the same, they are
algebraically independent over $\Q$.

To prove $i\pi,\log\pi,\log\log\pi,\ldots, \log_{[n]}\pi$ are
$\Q$-algebraically independent, we use Schanuel's Conjecture
again.

Without loss of generality, we may assume the statement true for
$$
i\pi,\log\pi,\log\log\pi,\ldots, \log_{[n-1]}\pi
$$
 (by induction).

As before, any nontrivial $\Q$-linear relation among the
$i\pi,\log\pi,\ldots, \log_{[n]}\pi$ can be thought as a nontrivial
$\Z$-linear combination (by clearing denominators) and then as an
algebraic relation among their exponentials (by applying $\exp(.)$
at both sides).

More precisely:
$$
i\pi q+\sum_{k=1}^n q_k \log_{[k]}\pi =0
$$
 with $q,q_k\in\Z$ leads
us to
$$
(-1)^q\prod_{k=1}^n(\log_{[k-1]}\pi)^{q_k}=1
$$
 which equals
$$
\prod_{k=0}^{n-1}(\log_{[k]}\pi)^{q_{k+1}}=(-1)^q
$$
Since we are assuming $i\pi,\log\pi,\log\log\pi,\ldots,
\log_{[n-1]}\pi$ are $\Q$-algebraically independent (also
$\pi,\log\pi,\log\log\pi,\ldots, \log_{[n-1]}\pi$) this last
algebraic relation must be the trivial one, i.e. $q_k=0\:  \; \hbox{ for all } \;
1\leq k\leq n$ and $q$ even (but this is no so important).
Returning to our linear relation we get $i\pi q =0$ meaning $q=0$.

Therefore $A=\{i\pi,\log\pi,\log\log\pi,\ldots, \log_{[n]}\pi\}$ are
linearly independent over $\Q$ and by Schanuel's Conjecture, the
transcendence degree of $\Q(A, \exp(A))$ should be at least $n+1$.

Since $\exp(A)$ is algebraic over $\Q(A)$, this means
$$
\trdeg_{\Q}\Q(i\pi,\log\pi,\log\log\pi,\ldots, \log_{[n]}\pi)\geq
n+1,
$$
 i.e. $i\pi,\log\pi,\log\log\pi,\ldots, \log_{[n]}\pi$ are
algebraically independent over $\Q$ (then over $\QQ$ and hence
over $E$). %end of proof
\end{proof}

\begin{corollary}
The numbers
 $e, e^e, e^{e^e}, \ldots$ are $L-$algebraically independent.
\end{corollary}

% or if you prefer
%\begin{corollary}
%The numbers
% $e, e^e, \ldots$ are algebraically independent over $L$.
%\end{corollary}

\begin{proof} As before, we only have to prove they are so over $\Q$.
Again, this follows by induction.

Name $\exp^{[n]}(1)=\exp(\exp^{[n-1]}(1))$ and $\exp^{[0]}(1)=1$.

Let us assume the $\{\exp^{[k]}(1) \}_{k=1}^n$ are algebraically
independent over $\Q$. Then the set
$$
A=\{1,e,e^e,\ldots,\exp^{[n]}(1) \}=\{\exp^{[k]}(1) \}_{k=0}^n
$$
 is
$\Q$-linearly independent and by Schanuel's Conjecture we should
have
$$
n+1 \leq \trdeg_{\Q}\Q(A,\exp(A))=\trdeg_{\Q}\Q(\exp(A))
$$
because $A$ is algebraic over $\Q(\exp(A))$.

But $\exp(A)=\{\exp^{[k]}(1) \}_{k=1}^{n+1}$ would be algebraically
independent over $\Q$. This finishes the inductive step. %end of proof
\end{proof}

\section*{Acknowledgements}

We would like to express our appreciation to Professor Waldschmidt
for coordinating this project during the evening sessions, and
Georges Racinet as well, for giving his support and guidance. We
would also like to thank the organizers of the Arizona Winter
School 2008 without whom none of these would have been possible.

\bibliographystyle{amsplain}

\vfill

Chuangxun Cheng (Northwestern U.) \href{cxcheng@math.northwestern.edu}{cxcheng@math.northwestern.edu}

Brian Dietel (Oregon State U.)  \href{dietelb@onid.orst.edu}{dietelb@onid.orst.edu}

Mathilde Herblot  (U. Rennes) \href{mathilde.herblot@univ-rennes1.fr}{mathilde.herblot@univ-rennes1.fr}

Jingjing Huang (PennState U.) \href{huang@math.psu.edu}{huang@math.psu.edu}

Holly Krieger (U. Illinois at Chicago) \href{hkrieger@uic.edu}{hkrieger@uic.edu}

Diego Marques (U. Federal do Ceara) \href{diego@mat.ufc.br}{diego@mat.ufc.br}

Jonathan Mason (U. Wisconsin) \href{mason@math.wisc.edu}{mason@math.wisc.edu}

Martin Mereb  (U. Texas at Austin) \href{mmereb@gmail.com}{mmereb@gmail.com}

S. Robert Wilson (Santa Cruz, CA) \href{robwilson93@mac.com}{robwilson93@mac.com}

\end{document}